\documentclass[11pt,twoside]{article}
\usepackage{amsfonts}
\usepackage{amsbsy}
\usepackage{amssymb}
\usepackage{amsmath}
\usepackage{amsthm}

\usepackage{graphicx}

\usepackage{color}
\usepackage{verbatim}

\definecolor{light-gray}{gray}{0.55}

\newtheorem{definition}{Definition}
\newtheorem{theorem}{Theorem}
\newtheorem{lemma}{Lemma}
\newtheorem{corollary}{Corollary}
\newtheorem{remark}{Remark}
\newtheorem{proposition}{Proposition}

\begin{document}

\par\bigskip

\begin{center}
{\LARGE{{\bf Geometric series of positive linear operators and inverse Voronovskaya theorem}}}
\end{center}
\medskip

\centerline{\Large{\bf Ulrich Abel, Mircea Ivan and Radu P\u alt\u anea}} 

\bigskip

\bigskip

\begin{abstract}
{We define the associated geometric series for a large
class of positive linear operators and study the convergence of the series
in the case of sequences of admissible operators. We obtain
an inverse Voronovskaya theorem and we apply our results to the Bernstein operators, the Bernstein-Durrmeyer-type operators, and  the symmetrical version of Meyer-K\"onig and Zeller operators.}

\textit{Key words:} positive linear operators, iterates of operators, series of operators, antiderivative, inverse Voronovskaya operators, Bernstein operators, Bernstein-Durrmeyer type operators, symmetrical version of Meyer-K\" onig and Zeller operators.
\end{abstract}

\section{Introduction. Basic notions}

\par Let $X$ be a linear subspace of $C[0,1]$ and let $L:X\to X$ be a positive linear operator. Set $\mathbb{N}_0=\mathbb{N}\cup\{0\}$. Denote by $L^k$, $k\in\mathbb{N}_0$, the iterates of $L$, defined by $L^0=I$, where $I$ is the identical operator and $L^k=L\circ{\cdots}\circ L$, where $L$ appears $k$ times.  
\par While the iterates of $L$ are always well defined, the geometric series of $L$, namely
\begin{equation}\label{geom}
G_L=\sum_{k=0}^{\infty}L^k
\end{equation}
needs some restrictions. So, if there is $f\in X$, $f\not=0$, such that $L(f)=f$, then $G_L(f)$ is not defined.
\par There exists a rich bibliography concerning the convergence of iterates of positive linear operators, starting with the paper of 
Kelisky and Rivlin~\cite{KR} about the iterates of Bernstein operators. 
A {non-exhaustive} list of contributions in this direction is given in our references: \cite{AI}, \cite{ABC}, \cite{CS0}, \cite{CS}, \cite{GI}, \cite{GaIv1}, \cite{GR}, \cite{K}, \cite{GHR}, \cite{N}, \cite{NRW}, \cite{R1}, 
\cite{W}.

{Comparatively, the geometric series of positive linear operators have been rather neglected.}

In paper P\u alt\u anea~\cite{PR04} is made the first study of the convergence of the sequence of geometric series $G_{B_n}$, $n\in\mathbb{N}$, where $B_n$, $n\in\mathbb{N}$ are the Bernstein operators. 
{A different approach} is given in \cite{AIP}. A varint of convergence of geometric  series of Bernstein operators on a simplex was studied by Ra\c sa~\cite{Ra}.

{The convergence of geometric series attached to Bernstein-Durrmeyer operators was  studied by Abel in \cite{Ab09}. 
\par In the present paper we consider geometric series
{related} to a large class of positive linear operators 
{including the majority of the usual approximation operators which preserve linear functions}. 
}

{Throughout the paper we will use the following notations:
\begin{description}\setlength{\itemsep}{0ex}
\setlength{\leftmargini}{3em}
\setlength{\itemindent}{-2em}
\setlength{\labelsep}{0em}
\item[{\makebox[4em][l]{$B[0,1]$,}}] 
the Banach space of all bounded functions $f\colon [0,1]\to\mathbb{R}$, endowed with the sup~norm $\|f\|=\sup_{x\in[0,1]}|f(x)|$;
\item[{\makebox[4em][l]{$C[0,1]$,}}]the Banach space of all continuous functions $f\colon [0,1]\to\mathbb{R}$, endowed with the same sup~norm;
\item[{\makebox[4em][l]{$C^k[0,1]$,}}]the space of all functions  $f\colon [0,1]\to\mathbb{R}$, possessing a continuous derivative of order $k$, $k\in \mathbb{N}_0$, where $C^0[0,1]=C[0,1]$;
\item[{\makebox[4em][l]{$C^k(0,1)$,}}] the set of all functions $f:[0,1]\to\mathbb{R}$ which have a continuous $k$-th derivative on the interval $(0,1)$;
\item[{\makebox[4em][l]{$\psi$,}}] the function $\psi(x)=x(1-x)$, $x\in[0,1]$;
\item[{\makebox[4em][l]{$\Pi_k$,}}]the set of all polynomial functions of degree at most $k$ defined on $[0,1]$;
\item[{\makebox[4em][l]{$e_j$,}}]the monomials $e_j(t)=t^j$ on $[0,1]$, $j=0,1,2,\ldots$.
\end{description}
}

\par Consider the following linear space of functions:
$$C_{\psi}[0,1]:=\{{f}:[0,1]\to\mathbb{R}\,|\;\exists {g}\in B[0,1]\cap C(0,1):\; f=\psi g\}.$$
Equivalently, we can write
$$C_{\psi}[0,1]:=\{f\in C[0,1]\,|\;\exists M>0:\;|f(x)|\le M\psi(x),\;x\in[0,1]\}.$$

Note that $\psi C[0,1]\subset C_{\psi}[0,1]$. On space $C_{\psi}[0,1]$ we consider the norm 
$$\|f\|_{\psi}:=\sup_{x\in(0,1)}\frac{|f(x)|}{\psi(x)},\;f\in C_{\psi}[0,1].$$
For simplicity,  we write sometimes $C_{\psi}[0,1]$ instead of $(C_{\psi}[0,1],
\|\cdot\|_{\psi})$. 
The following proposition can be immediately proved.
\begin{proposition}
The linear space $C_{\psi}[0,1]$ endowed with the norm $\|\cdot\|_{\psi}$ is a Banach space. 
\end{proposition}  
\begin{remark}\label{R1}
{\rm However, the linear space $C_{\psi}[0,1]$ endowed with the norm 
$\|\cdot\|$ is not a Banach space. For instance the sequence of functions $(\psi f_n)_n$, where $f_n$ are defined by $f_n(x)=n\sqrt{n}\cdot x$, for $x\in [0,\frac1n]$ and $f_n(x)=\frac1{\sqrt{x}}$, for $x\in[\frac1n,1]$ is fundamental, but converges uniformly to $\sqrt{x}(1-x)$ which is not an element of $C_\psi[0,1]$.} \end{remark} 

\par One of the basic properties of linear positive operators which preserve linear functions is given in the following well-known {result}. 

\par {\bf Theorem A}~[see, e.g., \cite[p.~5]{PR04}] {\it  If $L:C[0,1]\to C[0,1]$ is a positive linear operator such that $L(e_j)=e_j$, for $j=0,1$, then we have 
$$L(f)(0)=f(0) \ \mbox{and} \ L(f)(1)=f(1), \ \mbox{for all} \ f\in C[0,1].$$
}
\par If $X$ and $Y$ are two subsets of $C[0,1]$, denote $X+Y=\{f+g|\,f\in X,g\in Y\}.$ 
Denote by $B_1$, the Bernstein operator of order one: 
\begin{equation}\label{ee3}
B_1(f)(x)=(1-x)f(0)+x f(1),\;f\in C[0,1],\; x\in[0,1].
\end{equation}
\par From the well-known theorem of 
Korovkin~\cite[Theorem 8, Chapter 1, \S5]{Ko} we can derive 
the following immediate consequence.

\par\ \\
{\bf Theorem B.} {\it If $L:C[0,1]\to C[0,1]$ is a positive linear operator such that $L(\varphi_j)=\varphi_j$, $j=1,2,3$, where $\{\varphi_1, \varphi_2, \varphi_3\}$ is a Chebyshev system on interval $[0,1]$, then $L=I$.} 

\section{Operators from the class $\Lambda$ and their iterates}

Let $L:C[0,1]\to C[0,1]$ be a positive linear operator  preserving constants.
By using the Jessen inequality for functionals~\cite{Jessen-31}, we obtain
\begin{equation}\label{Eq.Lpsi<=psi}
0\le L(\psi) \le \psi.
\end{equation}
Eq.~(\ref{Eq.Lpsi<=psi}) yields

\begin{equation}\label{Eq.||Lpsi||<=1}
{\|L(\psi)\|_{\psi}}\le 1.
\end{equation} 
The previous equations imply that the function 
$${b_L}(x)=\left\{\begin{array}{cl}
\frac{L(\psi)(x)}{\psi(x)},& x\in(0,1),\\
{\|L(\psi)\|_{\psi}},&x\in\{0,1\},\end{array}\right.$$
belongs to $B[0,1]\cap C(0,1)$ and 
\begin{equation}\label{Eq.Lpsile|Lpsi|psi}
L(\psi) = {b_L}\psi\le{\|L(\psi)\|_{\psi}}\psi.
\end{equation}

\begin{definition}\label{D1}
Denote by $\Lambda$ the class of all positive linear operators 
$L:C[0,1]\to C[0,1]$ satisfying the following conditions:
\begin{itemize}
\item[{$(a)$}] $L$ preserves linear functions;
\item[{$(b)$}] ${\|L(\psi)\|_{\psi}}<1$;
\item[{$(c)$}] $L\not=B_1$.
\end{itemize}
\end{definition}

\begin{lemma}\label{L0}
If $L\in\Lambda$, then the following properties are satisfied:
\begin{itemize}
\item[$(i)$]$\|L(\psi)\|_{\psi}>0$; 
\item[$(ii)$] $L(C_{\psi}[0,1])\subset C_{\psi}[0,1]$;
\item[$(iii)$] 
$\|L\|_{\mathcal{L}(C_{\psi}[0,1],C_{\psi}[0,1])}\;
= \;\|L(\psi)\|_{\psi}=\|b_L\|$.
\end{itemize}

\end{lemma}
\par{\bf Proof.}

$(i)$ Suppose, for reductio ad absurdum, that  $\|L(\psi)\|_{\psi}=0$. It follows that $b_L=0$, hence $L(\psi)=0$. We show 
that this implies  $L=B_1$. Indeed, let $f\in C^1[0,1]$ and denote $M=2\|(f-B_1(f))'\|$. Since $f(0)=B_1(f)(0)$ we obtain for  $t\in [0,\frac12]$ that
$|f(t)-B_1(f)(t)|\le \|(f-B_1(f))'\|t\le M\psi(t)$.  The same inequality can be proved for $t\in [\frac12,1]$. Then, for $t\in[0,1]$, we deduce $|L(f)(t)-B_1(f)(t)|=|L(f-B_1(f))(t)|\le L(|f-B_1(f)|)(t)\le M L(\psi)(t)=0.$  Hence $L(f)=B_1(f)$, for $f\in C^1[0,1]$. Since $C^1[0,1]$ is dense in $C[0,1]$ it follows $L=B_1$. But this contradicts the hypothesis. It follows that $\|L(\psi)\|_{\psi}>0$.

\par $(ii)$ Let $g\in B[0,1]\cap C(0,1)$. 
We have $|L(\psi g)|\le \|g\|L(\psi)\le \|g\|\, \|L(\psi)\|_{\psi}\,\psi$.
It follows that $L(\psi g)\in C_{\psi}[0,1]$. 
\par $(iii)$ From above, for $g\in B[0,1]\cap C(0,1)$, we obtain:
$$\sup_{x\in(0,1)}\frac {|L(\psi g)(x)|}{\psi(x)}
\le \|L(\psi)\|_{\psi} 
\sup_{x\in(0,1)}|g(x)| =\|L(\psi)\|_{\psi}\;\|\psi g \|_{\psi},$$
i.e.,
$$\|L(f)\|_\psi \le \|L(\psi)\|_{\psi} \; \|f\|_\psi, \qquad f\in C_\psi[0,1]. $$
But $\|L(\psi)\|_{\psi}=\|b_L\|$. For $f=\psi$,  we have equality and the proof is complete.
\hfill $\Box$

\begin{definition}\label{D2}
Denote by $\Lambda_0$, the class of positive linear operators $L:C[0,1]\to C[0,1]$ which satisfy the following conditions:
\begin{itemize}
\item[$(a_1)$] $L$ preserves linear functions;
\item[$(b_1)$] $L(\Pi_2)\subset \Pi_2$;
\item[$(c_1)$] $L\not= I\;\mbox{and}\;L\not= B_1$.
\end{itemize}
\end{definition} 

\begin{lemma}\label{L1}
We have $\Lambda_0\subset\Lambda$. Moreover, $\Lambda_0$ coincides with the class of those operators $L\in\Lambda$, for which ${b_L}$ 
is a constant function with value in interval $(0,1)$.
\end{lemma}

\par{\bf Proof.} Let $L\in\Lambda_0$. By (\ref{Eq.Lpsile|Lpsi|psi}), we have $L(\psi)=b_L\,\psi$.
Since $L\psi\in \Pi_2$, we deduce that there exist constants
$a$, $b$, $c$ such that $b_L\psi= a+b\,e_1+ c\,\psi$. From Theorem A it follows $b_L(x)= c,$ $x\in(0,1)$.
Taking into account the definition of $b_L$, we deduce that
$b_L=c\, e_0$.

From relation (\ref{Eq.Lpsi<=psi}) we obtain $c\le 1$. If $c=1$, then operator $L$ preserves functions $e_0,e_1,\psi$ which form a Chebyshev system on interval $[0,1]$ and by Theorem B it follows $L=I$. Contradiction to $(c_1)$. Therefore ${\|L(\psi)\|_{\psi}}=c<1$. If $c=0$, which is equivalent to condition $\Vert L(\psi)\Vert_{\psi}=0$, then like in the proof of Lemma \ref{L0} we obtain $L=B_1$. Contradiction. Consequenly we have $c\in(0,1)$ and $L\in\Lambda$. 
\par Conversely, if for $L\in\Lambda$, there is a constant $c\in(0,1)$ such that ${b_L}(x)=c$, for $x\in[0,1]$, it follows that $L$ satisfies condition $(b_1)$, since the functions $e_0,e_1,\psi$ forms a basis of $\Pi_2$. Also, since $L(\psi)\not=\psi$ we have $L\not= I$. \hfill
$\Box$

\begin{theorem}\label{TA1}
If $L\in\Lambda$ we have
\begin{equation}\label{eA1}
\lim_{k\to\infty}\|L^k(f)-B_1(f)\|_{\psi}=0,\qquad \text{ for all } f\in C_{\psi}[0,1]+\Pi_1.
\end{equation}
\end{theorem}
\par{\bf Proof.}
Let $f=h+P_1$,  where $h\in C_\psi[0,1]$ and $P_1\in \Pi_1$.
Since $L^k$ and $B_1$ preserve $P_1$ we have:
$$L^k(f)- B_1(f)= L^k(h)- B_1(h) = L^k( h- B_1(h)).$$
Since $h- B_1(h)$ belongs to $C_\psi[0,1]$, 
by using Lemma~\ref{L0}~$(ii)$ and Eq.~(\ref{Eq.Lpsile|Lpsi|psi}) repetitively, we obtain:
$$\|L^k(f)- B_1(f)\|_\psi = \|L^k( h- B_1(h))\|_\psi
\le \big(\|L(\psi)\|_\psi\big)^k \;\|h-B_1h\|_\psi.$$
 Since $\|L(\psi)\|_\psi<1$, passing to limit as $k\to\infty$ the proof of (\ref{eA1}) is complete. \hfill $\Box$

\par We mention a different result for the convergence of iterates of operators given by Gavrea and Ivan \cite{GI}, which in a slightly modified form could be read as follows. 
\par\ \\
{\bf Theorem C} {\it If $L:C[0,1]\to C[0,1]$ is a positive linear operator such that $L(e_j)=e_j$, $j=0,1$, and there is $m>0$ and $q\ge1$ such that $\psi-L(\psi)\ge m\,\psi^q$, then we have 
\begin{equation}\label{eGI1}
\lim_{k\to\infty}\|L^k(f)-B_1(f)\|=0,\;\mbox{for all}\;f\in C[0,1].
\end{equation}}
Note that, for $q=1$, the conditions in Theorem C are satisfied by operators in class $\Lambda$. 

However Theorem \ref{TA1} is  not a consequence of Theorem C, since it is given with respect to norm $\|\cdot\|_{\psi}$.

\section{Geometric series of operators from class $\Lambda$}

\par We consider the geometric series $G_L$ attached to an operator $L\in\Lambda$ on the space $C_{\psi}[0, 1]$.

\begin{theorem}\label{T2}
For any $L\in\Lambda$ we have:
\begin{itemize}
\item[$(i)$] the operator $G_L:C_{\psi}[0,1]\to C_{\psi}[0,1]$, given in (\ref{geom}) is well defined if we consider the convergence with regard to the norm 
$\|\cdot\|_{\psi}$ and hence also with regard to norm $\|\cdot\|$;
\item[$(ii)$] the operator $G_L$ is positive and linear;
\item[$(iii)$] $(1-\Vert b_L\Vert )G_L(\psi)\le\psi$;
\item[$(iv)$] $\|G_L\|_{{\cal L}(C_{\psi}[0,1],C_{\psi}[0,1])}\le
({1-\|b_L\|})^{-1}$.
\end{itemize}
\end{theorem}

 \par {\bf Proof.} $(i)$ Let $f\in C_{\psi}[0,1]$. From Lemma~\ref{L0}~$(ii)$ it follows by induction that $L^k(f)\in C_{\psi}[0,1]$, for all $k\ge 0$. Hence the partial sums of the series $\sum_{k=0}^{\infty}L^k(f)$ are in $C_{\psi}[0,1]$. From Lemma~\ref{L0}~$(iii)$ we obtain $\|L^k\|_{{\cal L}(C_{\psi}[0,1],C_{\psi}[0,1])}
\le \|b_L\|^k$, for $k\ge0$. It follows that
\begin{equation}\label{Eq.Th2.(i)}
\sum_{k=0}^{\infty}\|L^k\|_{{\cal L}(C_{\psi}[0,1],C_{\psi}[0,1])}\le 
\frac{1}{1-\|b_L\|}
\end{equation}%
(note that $\|b_L\|<1$  by Definition~\ref{D1}~$(b)$).
Consequently, the series $\sum_{k=0}^{\infty}L^k(f)$ is convergent in the space $C_{\psi}[0,1]$, with regard to the norm $\|\cdot\|_{\psi}$. 
\par $(ii)$ It is immediate.
\par $(iii)$ It follows from Eq.~(\ref{Eq.Lpsile|Lpsi|psi}).  
\par $(iv)$ The inequality is a consequence of Eq.~(\ref{Eq.Th2.(i)}). \hfill $\Box$

\begin{theorem}\label{T3}
For any operator $L\in\Lambda$, the following equalities are true on the Banach space $(C_{\psi}[0,1],\|\cdot\|_{\psi})$:
\begin{eqnarray}
(I-L)\circ G_L&=&I;\label{e7}\\
G_L\circ(I-L)&=&I.\label{e8}
\end{eqnarray}
\end{theorem}

\par {\bf Proof.} Since $\|L^k\|_{{\cal L}(C_{\psi}[0,1],C_{\psi}[0,1])}\le \|b_L\|^k$, for $k\in\mathbb{N}_0$,  we obtain 
$$\lim_{k\to\infty}\|L^k(f)\|_{\psi}=0,\qquad\text{ for all }f\in C_{\psi}[0,1].$$%
Then, using the continuity of operator $I-L$ we obtain successively, for  $f\in C_{\psi}[0,1])$:
\begin{eqnarray}
((I-L)\circ G_L)(f)&=&(I-L)\left(\lim_{m\to\infty}\sum_{k=0}^m L^k(f)\right)\nonumber\\
&=&\lim_{m\to\infty}(I-L)\left(\sum_{k=0}^m L^k(f)\right)\nonumber\\
&=&\lim_{m\to\infty}(I-L^{m+1})(f)\nonumber\\
&=&f.\nonumber
\end{eqnarray}
For $f\in C_{\psi}[0,1]$, we have
$$(G_L\circ(I-L))(f)=\lim_{m\to\infty}\sum_{k=0}^m (L^k\circ (I-L))(f)= \lim_{m\to\infty}(I-L^{m+1})(f)=f.\quad\Box$$

\section{Convergence of geometric series on the space $C_{\psi}[0,1]$}

\par For $f\in C(0,1)\cap B[0,1]$ and $x\in[0,1]$ define
\begin{equation}\label{e6}
F(f)(x)=(1-x)\int_0^x tf(t)dt+x\int_x^1(1-t)f(t)dt.
\end{equation}

\par In \cite{PR06}, the following immediate result is proved.
\par\ \\
{\bf Lemma E} {\it For any $f\in C[0,1]$, we have $F(f)\in \psi C[0,1]\cap C^2[0,1]$ and
\begin{equation}
F''(f)=-f.
\end{equation}}

A slightly  modified version of it is the following.

\begin{lemma}\label{L4}
For any $f\in B[0,1]\cap C(0,1)$, we have $F(f)\in \psi C[0,1]\cap C^2(0,1)$ and
\begin{equation}\label{e9a}
F''(f)(x)=-f(x),\qquad x\in(0,1).
\end{equation}
\end{lemma}
\par {\bf Proof.} Let $f\in B[0,1]\cap C(0,1)$. The function $g:=F(f)\psi^{-1}$ defined and continuous on $(0,1)$  can be extended by continuity  to the  interval $[0,1]$, since $\lim\limits_{x\to 0+}g(x)=\int_0^1(1-t)f(t)dt$ and $\lim\limits_{x\to 1-}g(x)=\int_0^1tf(t)dt$. 
 Relation (\ref{e9a}) follows immediately. \hfill $\Box$

\par Consider a sequence of operators $(L_n)_n$, $L_n\in \Lambda$, $n\in\mathbb{N}$. For simplicity, we denote 
\begin{eqnarray}
G_n&:=&G_{L_n},\;\alpha_n:=1-b_{L_n},\;\nu_n:=\min_{x\in[0,1]}\alpha_n,\label{e40}\\ \eta_n&:=&\sup_{x,y\in[0,1]}|\alpha_n(x)-\alpha_n(y)|/\nu_n,\label{e41}\\
M^k_n(x)&:=& L_n((e_1-xe_0)^k,x),
\end{eqnarray}
for $x\in[0,1]$, $n\in\mathbb{N}$, $k\in\mathbb{N}_0$. Then $M_n^2=\alpha_n\psi$. With these notations we have:

\begin{theorem}\label{T4}
Let the sequence of operators $(L_n)_n$, $L_n\in\Lambda$, $n\in\mathbb{N}$, which satisfies the following ``little-o''  conditions:
\begin{equation}\label{e22}
M^4_n(x)={\rm o}(M_n^2(x)),\;\mbox{uniformly with regard to}\;x\in[0,1]; 
\end{equation}
\begin{equation}\label{e21}
\eta_n={\rm o}(1),
\end{equation}
as $n\to\infty$.
Then, for any $f\in B[0,1]\cap C(0,1)$, we have
\begin{equation}\label{e23}
\lim_{n\to\infty}\|\alpha_n G_n(\psi f)-2F(f)\|_{\psi}=0.
\end{equation}\end{theorem}

\par{\bf Proof.} We split the proof in two steps.

\par {\bf Step 1.} We prove relation (\ref{e23}) in the case $f\in C[0,1]$. Set $g=\psi f$.
Fix $n\in\mathbb{N}$. Denote for simplicity $f_n:=f/\alpha_n$ and $H_n:=F(f_n)$. It follows $f_n\in B[0,1]\cap C(0,1)$ and using Lemma \ref{L4} we have
$H_n\in C^2(0,1)\cap C_{\psi}[0,1]$. We obtain for  $y\in(0,1)$ and $t\in[0,1]$:
$$H_n(t)=H_n(y)+H_n'(y)(t-y)+\frac12\cdot H_n''(y)(t-y)^2+\Theta_{n,y}(t),$$
where
$$\Theta_{n,y}(t):=\int_y^t(t-u)(H_n''(u)-H_n''(y))\,du,$$
and $H_n''(t)=-f_n(t)$, for $t\in (0,1)$.
Fix $y\in(0,1)$.  Applying the linear functional $L_n(\cdot,y)$, and taking into account condition~$(a)$ from Definition \ref{D1} we obtain 
$$L_n(H_n)(y)=H_n(y)-\frac12\, f_n(y) M_n^2(y)+L_n(\Theta_{n,y})(y).$$
It follows
\begin{eqnarray}
(I-L_n)(H_n)(y)&=&\frac12\, g(y)-L_n(\Theta_{n,y})(y).\label{e12}
\end{eqnarray}
Note that Eq.~(\ref{e12}) which was derived for $y\in (0,1)$ is also true for $y=0$ or $y=1$, since in these cases the both side terms of the equation vanish.
\par Denote $S_n(y)=L_n(\Theta_{n,y})(y)$, $y\in[0,1]$.  Since $H_n\in C_{\psi}[0,1]$, we obtain $(I-L_n)(H_n)\in C_{\psi}[0,1]$ and then we deduce $S_n\in C_{\psi}[0,1]$. 
Applying operator $G_n$, to the functions in both sides of Eq.~(\ref{e12}) and taking into account Theorem~\ref{T3} we obtain, for $x\in[0,1]$:
$$ H_n(x)=\frac 12\, G_n(g)(x)- G_n(S_n)(x).$$
Hence
\begin{equation}\label{e15}
\|\alpha_n\big(G_n(g)-2F(f_n)\big)\|_{\psi}=2\|\alpha_nG_n(S_n)
\|_{\psi}.
\end{equation}
\par Let $\omega (f,\delta):=\sup\{|f(u)-f(v)|,\;u,v\in [0,1],\;|u-v|\le\delta\}$, for $\delta>0$, be the modulus of continuity of $f$. 
By condition~(\ref{e22}), there is a sequence $(\beta_n)_n$ of positive numbers such that $\beta_n\to 0$ $(n\to\infty)$ and $M_n^4(x)\le\beta_n M_n^2(x)$, for all $x\in[0,1]$ and $n\in\mathbb{N}$.

For $u,y\in[0,1]$,  we obtain
\begin{eqnarray}
|f_n(u)-f_n(y)|&\le&\Big|f(u)\Big(\frac1{\alpha_n(u)}-\frac1{\alpha_n(y)}\Big)\Big|+\frac{|f(u)-f(y)|}{\alpha_n(y)}\nonumber\\
&\le&\frac{\eta_n\|f\|}{\alpha_n(y)}+\frac1{\alpha_n(y)}\Big(1+\frac{(u-y)^2}{\beta_n}\Big)\;\omega(f,\sqrt{\beta_n}).\nonumber
\end{eqnarray}
Consequently, for $t\in[0,1]$, by Lemma~\ref{L4}, we have: 
\begin{eqnarray}
|\Theta_{n,y}(t)|&\le&\Big|\int_y^t(t-u)|f_n(u)-f_n(y)|\,du\Big|\nonumber\\
&\le&\frac1{2\alpha_n(y)}\Big[\eta_n\|f\|(t-y)^2+\Big((t-y)^2+\frac{(t-y)^4}{6\beta_n}\Big)\omega(f,\sqrt{\beta_n})\Big].\nonumber
\end{eqnarray}
Then we deduce:
\begin{eqnarray}
|S_n(y)|=|L_n(\Theta_{n,y})(y)|&\le&\frac1{2\alpha_n(y)}\Big[\eta_n\|f\|M_n^2(y)+\Big(M_n^2(y)+\frac1{6\beta_n}\cdot M_n^4(y)\Big)\omega(f,\sqrt{\beta_n})\Big]\nonumber\\
&\le&\psi(y)\Big[\frac12\eta_n\|f\|+\frac7{12}\omega(f,\sqrt{\beta_n})\Big],\nonumber
\end{eqnarray}
because $M_n^2=\alpha_n\psi$.  Therefore, by Th.~\ref{T2}~$(iii)$, it follows that 
\begin{eqnarray}
\|\alpha_n G_n(S_n)\|_{\psi}&\le&\left(\frac12\eta_n\|f\|+\frac7{12}\omega(f,\sqrt{\beta_n})\right)
\|\alpha_n G_n(\psi)\|_{\psi}\nonumber\\
&\le&\left(\frac12\eta_n\|f\|+\frac7{12}\omega(f,\sqrt{\beta_n})\right) \frac{\Vert\alpha_n\Vert}{\nu_n}\nonumber\\
&\le&\left(\frac12\eta_n\|f\|+\frac7{12}\omega(f,\sqrt{\beta_n})\right)(1+\eta_n).\label{ee24}
\end{eqnarray}

Since $f\in C[0,1]$ we have $\lim_{\rho\to0}\omega(f,\rho)=0$. From relations (\ref{e22}), (\ref{e21}), (\ref{e15}) and (\ref{ee24}) it results:
\begin{equation}\label{ee23}
\lim_{n\to\infty}\|\alpha_n \big(G_n(g)-2F(f_n)\big)\|_{\psi}=0.
\end{equation}
\par Finally, we use the inequality
\begin{eqnarray}
\|\alpha_n G_n(g)-2F(f)\|_{\psi}&\le&\|\alpha_n \big(G_n(g)-2F(f_n)\big)
\|_{\psi}\nonumber\\
&&+2\|\alpha_nF(f_n)-F(f)\|_{\psi}.\label{e25}
\end{eqnarray}
For $y\in(0,1),$ we obtain
\begin{eqnarray}
|\psi^{-1}(y)\big(\alpha_n(y)F(f_n)(y)-F(f)(y)\big)|&=&\Big|\psi^{-1}(y)F\Big(f\cdot\frac{\alpha_n(y)e_0-\alpha_n}{\alpha_n}\Big)(y)\Big|\nonumber\\
&\le& \psi^{-1}(y)F\Big(|f|\cdot\frac{|\alpha_n(y)e_0-\alpha_n|}{\alpha_n}\Big)(y)\nonumber\\
&\le&\eta_n\psi^{-1}(y)F(|f|)(y)\nonumber\\
&\le&\eta_n\|F(|f|)\|_{\psi}.\nonumber
\end{eqnarray}
Hence, by using assumption (\ref{e21}) it follows
that\begin{equation}\label{e26}
\lim_{n\to\infty}\|\alpha_nF(f_n)-F(f)\|_{\psi}=0.
\end{equation}
From Eqs.~(\ref{e25}), (\ref{ee23}) and (\ref{e26}) we obtain (\ref{e23}).

\par{\bf Step 2.} Now we prove relation (\ref{e23}) in the general case when $f\in B[0,1]\cap C(0,1)$. Let  $\varepsilon$, $0<\varepsilon<1$. Then choose a number  $\delta$ with $0<\delta<\frac{\varepsilon}{96(\|f\|+1)}$. Define

$$f_{\delta}(t):=\left\{\begin{array}{ll}f(\delta),&t\in[0,\delta]\\ f(t),&t\in[\delta,1-\delta]\\f(1-\delta),&t\in[1-\delta,1].\end{array}\right.,\qquad
\varphi_{\delta}^1(t):=\left\{\begin{array}{ll}1,&t\in[0,\delta]\\2-t/\delta,&t\in[\delta,2\delta]\\0,&t\in[2\delta,1]\end{array}\right..$$
Then define the function $\varphi_{\delta}^2(t)=\varphi_{\delta}^1(1-t)$, $t\in [0,1]$. 
Note that:
$$0\le F(\varphi_{\delta}^1)(x)\le(1-x)\int_0^{\min\{x,2\delta\}}tdt+x\int_{\min\{x,2\delta\}}^{2\delta}(1-t)dt\le (1-x)x\delta+2x(1-x)\delta.$$ 
It follows $\|F(\varphi_{\delta}^1)\|_{\psi}\le 3\delta$. In an analogous mode we have $\|F(\varphi_{\delta}^2)\|_{\psi}\le 3\delta$. 

\par Since $f_{\delta},\varphi_{\delta}^1,\varphi_{\delta}^2\in C[0,1]$, using Step 1, there is $n_{\varepsilon}\in\mathbb{N}$, such that, for any
 integer $n\ge n_{\varepsilon}$, we have
$$\|\alpha_n G_n(\psi f_{\delta})-2F(f_{\delta})\|_{\psi}<\frac{\varepsilon}3,\;\mbox{and}\;
\|\alpha_n G_n(\psi \varphi_{\delta}^j)-2F(\varphi_{\delta}^j)\|_{\psi}<\frac{\varepsilon}{24(\|f\|+1)},\quad j=1,2.$$ 
\par Let now  $n\ge n_{\varepsilon}$. We have
\begin{eqnarray}
\|\alpha_n G_n(\psi f)-2F(f)\|_{\psi}&\le&\|\alpha_n G_n(\psi f_{\delta})-2F(f_{\delta})
\|_{\psi}\nonumber\\
&&+\|\alpha_nG_n(\psi(f-f_{\delta}))\|_{\psi}+2\|F(f-f_{\delta})\|_{\psi}\nonumber\\
&<&\frac{\varepsilon}3+\sum_{j=1}^2 \|\alpha_n G_n(\psi |f-f_{\delta}|\varphi_{\delta}^j)
\|_{\psi}+2\|F(f-f_{\delta})\|_{\psi}\nonumber\\
&\le&\frac{\varepsilon}3+2\|f\|\sum_{j=1}^2\|\alpha_n G_n(\psi\varphi_{\delta}^j)\|_{\psi}+4 \|f\|\sum_{j=1}^2\|F(\varphi_{\delta}^j)\|_{\psi}\nonumber\\
&\le&\frac{\varepsilon}3+2\|f\|\sum_{j=1}^2\Big[\|\alpha_n G_n(\psi\varphi_{\delta}^j)
-2F(\varphi_{\delta}^j)\|_{\psi}+4\|F(\varphi_{\delta}^j)\|_{\psi}\Big]\nonumber\\
&\le&\frac{\varepsilon}3+2\|f\|\Big[2\cdot\frac{\varepsilon}{24(\|
f\|+1)}+8\cdot 3\delta\Big]\nonumber\\
&<&\varepsilon.\nonumber
\end{eqnarray}
Since $\varepsilon>0$ was arbitrarily chosen it follows relation (\ref{e23}) in the general case when $f\in B[0,1]\cap C(0,1)$. \hfill 
 $\Box$

\begin{corollary}\label{C1}
In conditions of Theorem \ref{T4} there is a constant $M>0$, such that
\begin{equation}
\|\alpha_n G_n\|_{\mathcal{L}(C_{\psi}[0,1],C_{\psi}[0,1])}\le M,
\qquad \text{ for all } n\in\mathbb{N}.
\end{equation}
\end{corollary}
 
\par{\bf Proof.}  We can apply the 
uniform boundedness principle
to operators $\alpha_n G_n:C_{\psi}[0,1]\to C_{\psi}[0,1]$ by taking into account relation (\ref{e23}) and the fact that $C_{\psi}[0,1]$ is a Banach space. \hfill $\Box$

\begin{remark}
{\rm By taking into account Remark \ref{R1} we cannot derive by the 
uniform boundedness principle a bound for $\|\alpha_n G_n\|_{\mathcal{L}((C_{\psi}[0,1],
\|\cdot\|),(C_{\psi}[0,1],\|\cdot\|))}$, for all $n\in\mathbb{N}$.}
\end{remark} 

\section{Inverse Voronovskaya theorem}

\par An inverse Voronovskaya theorem was established using the theory of semigroups of operators, for instance in \cite{AD}. Here we give a variant of inverse Voronovskaya theorem using the geometric series of operators. We use the notations given in (\ref{e40}) and (\ref{e41}).
\par The Voronvskaya theorem can be expressed in the strong form of the convergence in Banach space $(C_{\psi}[0,1],\|\cdot\|_{\psi})$. From the general result given in \cite[Corollary 4.3]{GP}, we deduce, with the notations in this paper: 
\par\ 
\par {\bf Theorem D.} {\it Let a sequence of positive linear operators $(L_n)_n$, $L_n\in \Lambda_0$, with the additional conditions: $L(C[0,1])\subset C^4[0,1]$,  $L(\Pi_j)\subset\Pi_j$, $j=3,4$ and relation (\ref{e22}) is true. 
Then, for any $f\in C^2[0,1]$, we have:
\begin{equation}\label{ed2}
\lim_{n\to\infty}\left\|\frac1{\nu_n}(L_n(f)-f)-\frac 12 f''\psi\right\|_{\psi}=0.
\end{equation}}

\par In this section we state a converse result to Theorem D.

 \begin{theorem}\label{T5} 
Let $(L_n)_n$ be a sequence of operators $L_n\in \Lambda$, $n\in\mathbb{N}$, which satisfy conditions (\ref{e22}), (\ref{e21}) and also the following condition:
\begin{equation}\label{e55}
L_n(\psi|\alpha-\alpha(x)e_0|)(x)={\rm o}(\nu_n^{2}\psi(x)),\;\mbox{uniformly for}\;x\in[0,1].
\end{equation}
If for $f\in C_{\psi}[0,1]+\Pi_1$  there holds   
\begin{equation}\label{ed3}
 \lim_{n\to\infty}\left\|\frac1{\alpha_n}(L_n(f)-f)-g\psi\right\|_{\psi}=0,
\end{equation}
with a certain $g\in C(0,1)\cap B[0,1]$, then $f\in C^2(0,1)$ and $g(x)=\frac12 f''(x)$, for $x\in(0,1)$.  
\end{theorem}

\par{\bf Proof.} Let $f\in C_{\psi}[0,1]+\Pi_1$ and $g\in C(0,1)\cap B[0,1]$ satisfying the conditions in the theorem. Put $f_1= f-B_1(f)$ and $h_n=\frac1{\alpha_n}(L_n(f)-f)$, for $n\in\mathbb{N}$. We have $f_1\in C_{\psi}[0,1]$ and $h_n\in C_{\psi}[0,1]$. Then, there is $g_1\in C(0,1)\cup B[0,1]$, such that $f_1=\psi g_1$. First we prove the following limit 
\begin{equation}\label{ed7}
\lim_{n\to\infty}\|\alpha_n G_n(h_n)+f_1\|_{\psi}=0. 
\end{equation}
Indeed, we can write
$$\alpha_nG_n(h_n)=\alpha_n\sum_{k=0}^{\infty}L_n^{k}\left(\frac{L_n(f_1)-f_1}{\alpha_n}\right)=-f_1+\sum_{k=0}^{\infty}\alpha_nL_n^k\left(\frac1{\alpha_n}L_n(f_1)-L_n\left(\frac{f_1}{\alpha_n}\right)\right).$$

By condition~(\ref{e55}), there is a sequence $(\gamma_n)_n$ of positive numbers converging to zero such that 
$L_n(\psi|\alpha-\alpha(x)e_0|)(x)\le\gamma_n \nu_n^{2}\psi(x)$, for all $x\in[0,1]$ and $n\in\mathbb{N}$.

For $x\in[0,1]$ and $n\in\mathbb{N}$, we have:
\begin{eqnarray}
\left|\frac1{\alpha_n(x)}\cdot L_n(f_1)(x)-L_n\left(\frac{f_1}{\alpha_n}\right)(x)\right|&\le&L_n\left(|f_1|\cdot\frac{|\alpha_n-\alpha_n(x)e_0|}{\alpha_n\alpha_n(x)}\right)(x)\nonumber\\
&\le&\|g_1\|\gamma_n\psi(x).\nonumber
\end{eqnarray}
Then, using relation (iii) of Theorem \ref{T2}, we have:
\begin{eqnarray}
\left|\sum_{k=0}^{\infty}\alpha_nL_n^k\left(\frac1{\alpha_n}L_n(f_1)-L_n\left(\frac{f_1}{\alpha_n}\right)\right)\right|
&\le&\sum_{k=0}^{\infty}\alpha_n\|g_1\|\gamma_nL_n^k(\psi)\nonumber\\
&\le&\frac{\alpha_n}{\nu_n}\cdot\|g_1\|\gamma_n\psi\nonumber\\
&\le&(1+\eta_n)\gamma_n\|g_1\|\psi.\nonumber
\end{eqnarray}
We obtain (\ref{ed7}). Write
\begin{equation}\label{ed8}
\alpha_nG_n(h_n)=\alpha_nG_n(h_n-\psi g)+\alpha_nG_n(\psi g).
\end{equation}
By Theorem \ref{T4},  
\begin{equation}\label{ed9}
\lim_{n\to\infty}\|\alpha_nG_n(\psi g)-2 F(g)\|_{\psi}=0.
\end{equation}
\par From Corollary \ref{C1} there is a constant $M$ such that  $\|\alpha_nG_n
\|_{{\cal L}(C_{\psi}[0,1],C_{\psi}[0,1])}\le M$, $n\in\mathbb{N}$.  It follows $\|\alpha_nG_n(h_n-\psi g)\|_{\psi}\le M\|h_n-\psi g\|_{\psi}$ and from relation (\ref{ed3}) we obtain
\begin{equation}\label{ed10}
\lim_{n\to\infty}\|\alpha_nG_n(h_n-\psi g)\|_{\psi}=0.
\end{equation}
\par Now from relations (\ref{ed7}), (\ref{ed8}), (\ref{ed9}) and (\ref{ed10}) we obtain
$$-f_1=\lim_{n\to\infty}\alpha_nG_n(h_n)=2F(g).$$
Since $2F(g)\in C^2(0,1)$, it follows $-f_1\in C^2(0,1]$ and consequently $f\in C^2(0,1)$. Moreover, using Lemma \ref{L4} we deduce on interval $(0,1)$: $f''=(f_1+B_1(f))''=(f_1)''=(-2F(g))''=2g$. \hfill  $\Box$

\section{Applications}

\subsection{Bernstein operators}

The classical Bernstein operators are given by
$$B_n(f)(x)=\sum_{k=0}^{\infty}f\left(\frac kn\right)p_{n,k}(x),\;f\in C[0,1],\;x\in[0,1],$$  
where $p_{n,k}(x)=\binom{n}{k}x^k(1-x)^{n-k}$, $0\le k\le n$. 
\par If we denote the moments of operators $B_n$, by $M_n^k(x)$, $x\in[0,1]$, $k=0,1,\ldots$, then
$M_n^2(x)=\frac{\psi(x)}n,\;M_n^4(x)=\frac{\psi(x)}{n^2}\left[\left(3-\frac6n\right)\psi(x)+\frac1n\right].$
It follows $B_n\in\Lambda_0$, with $\alpha_n(x)=\nu_n=\frac1n$. Then, conditions (\ref{e21}) and (\ref{e55}) are automatically satisfied.  Also, for $n\ge2$, we have $M_n^4(x)/M_n^2(x)=\left(\frac 3n-\frac 6{n^2}\right)\psi(x)+\frac1{n^2}\le\frac 34\cdot\frac1n-\frac12\cdot\frac1{n^2}$. Hence relation (\ref{e22}) is satisfied too. Then we can apply Theorem \ref{T4} and Theorem \ref{T5} to Bernstein operators. 

\subsection{Bernstein-Durrmeyer type operators}

A class of Bernstein-Durrmeyer type operators $U_n^{\rho}$, $\rho\in(0,\infty)$, $n\in\mathbb{N}$, which preserve linear functions is considered in \cite{PR07}, \cite{PR10}: 
$$(U_n^{\rho}f)(x)=\sum_{k=1}^{n-1} F^{\rho}_{n,k}(f)p_{n,k}(x)+f(0)p_{n,0}(x)+f(1)p_{n,n}(x),\;f\in C[0,1],\;x\in[0,1]$$
with functionals $F^{\rho}_{n,k}$, $1\le k\le n-1,$ defined by: 
$$F_{n,k}^{\rho}(f):=\int_0^1 f(t)\mu^{\rho}_{n,k}(t)dt,\;
\mu^{\rho}_{n,k}(t) := \frac{t^{k\rho-1}(1-t)^{(n-k)\rho-1}}{B(k\rho,(n-k)\rho)},$$
where $B(x,y)$ is Euler's Beta function. 
\par We mention that in the special case $\rho=1$, one obtains the ``genuine'' Durrmeyer operator and for each $f\in C[0,1]$ we have $\lim_{\rho\to\infty}U^{\rho}_n(f)=B_n(f)$, $f\in C[0,1]$, uniformly. 
\par For fixed $\rho$, let denote the moments of operators $U_n^{\rho}$, by $M_n^k(x)$, $n\in\mathbb{N}$, $k=0,1,\ldots$, $x\in[0,1]$. In \cite{PR07}, \cite{PR10} the following formulas are derived:  
\begin{eqnarray}
&&M_n^0(x)=1,\nonumber\\
&&M_n^1(x)=0,\nonumber\\
&&M_n^2(x)=\frac{(\rho+1)x(1-x)}{n\rho+1},\nonumber\\
&&M_n^4(x)=\frac{3\rho(\rho+1)^2\psi^2(x)n}{(n\rho+1)(n\rho+2)(n\rho+3)}+\nonumber\\
&&+\frac{-6(\rho+1)(\rho^2+3\rho+3)\psi^2(x)+ (\rho+1)(\rho+2)(\rho+3)\psi(x)}
{(n\rho+1)(n\rho+2)(n\rho+3)}.\nonumber
\end{eqnarray} 
Hence $U_n^{\rho}\in\Lambda_0$ with $\alpha_n=\nu_n=\frac{\rho+1}{n\rho+1}$. Then relation (\ref{e21}) and (\ref{e55}) are automatically satisfied. Condition (\ref{e22}) is also satisfied, since 
$$M_n^4(x)/M_n^2(x)\le\frac{3\rho(\rho+1)n\psi(x)-6(\rho^2+3\rho+3)\psi(x)+(\rho+2)(\rho+3)}{(n\rho+2)(n\rho+3)}$$
and the functions given in the right side of this inequality converge uniformly to $0$. Consequently Theorem \ref{T4} and Theorem \ref{T5} can be applied to operators $U_n^{\rho}$.

\subsection{A symmetrical version of the Meyer-K\" onig and Zeller operators}

The Meyer-K\" onig and Zeller operators \cite{MKZ}, with a slight modification made by Cheney and Sharma \cite{CS} are defined, for $f\in C[0,1]$, by

\begin{equation}
Z_n(f)(x)= \left\{\begin{array}{ll}
\displaystyle\sum_{k=0}^{\infty}\binom{n+k}{k}(1-x)^{n+1}x^k f\left(\frac k{n+k}\right),&x\in[0,1),\\[4ex]
f(1),& x=1.\end{array}\right.
\end{equation}
These operators reproduce linear functions. 
\par We consider here a symmetric variant  of operators $Z_n$. Let $\tau:=e_0-e_1$ and define 
\begin{eqnarray}
Z_n^1(f)(x)&:=&Z_n(f\circ\tau)(1-x),\;f\in C[0,1],\quad x\in[0,1],\\
Z_n^{\star}(f)&:=&\frac12(Z_n+Z_n^1).
\end{eqnarray}

\par We have $Z_n^1(e_0)=Z_n(e_0)=e_0$ and, for $x\in[0,1]$: $Z_n^1(e_1)(x)=Z_n(e_0-e_1)(1-x)= 1-(1-x)=e_1(x)$. Then $Z_n^1$ and consequently also $Z^*_n$ preserve linear functions. 
\par Let $n\in\mathbb{N}$,  $k\in\mathbb{N}$ and $x\in[0,1]$. Denote $m_n^k(x)=Z_n((e_1-xe_0)^k)(x)$ and  $M_n^k(x)=Z^{\star}_n((e_1-xe_0)^k)(x)$. We have $Z_n^1((e_1-xe_0)^k)(x)=(-1)^km_n^k(1-x)$. Therefore, if $k$ is even we have
$M_n^k(x)=\frac12(m_n^k(x)+m_n^k(1-x)).$ 
Consider also the notations given in (\ref{e40}) and (\ref{e41}) for operators $L_n=Z_n^{\star}$, $n\in\mathbb{N}$. 
\par Becker and Nessel~\cite{BN} proved the double inequality: 
\begin{equation}\label{e67}
\frac{x(1-x)^2}{n+1}\left(1+\frac{2x}{n+2}\right)\le m_n^2(x)\le \frac{x(1-x)^2}{n+1}\left(1+\frac{2x}{n+1}\right).
\end{equation}
From this we deduce
$$\frac{\psi}{2(n+1)}\left(1+\frac{4\psi}{n+2}\right)\le M_n^2(x)\le \frac{\psi}{2(n+1)}\left(1+\frac{4\psi}{n+1}\right).$$
Hence $M_n^2(x)=\alpha_n\psi$, $\alpha_n\in B[0,1]\cap C(0,1)$ and $\frac1{2(n+1)}\le \alpha_n(x)\le \frac{n+2}{2(n+1)^2}$, for $x\in[0,1]$. It follows $Z_n^{\star}\in\Lambda.$ Moreover, it follows that $\eta_n\le\frac1{n+1}$, hence, the sequence of operators $(Z_n^{\star})_n$ satisfies condition (\ref{e21}).

\par We start now to estimate the values of $Z_n(e_r)$, for an arbitrary $r\ge 2$. A complete asymptotic expansion is given by Abel in \cite{Ab95}. Here we need a more simple estimate, but in which the factor $\psi$ is taken into account. In particular case we obtain a new estimate of moment $m_n^2(x)$, which, from certain point of view is more precise then the estimate of Becker and Nessel (\ref{e67}). 

\par Let $x^{\underline{j}}=x(x-1)\ldots (x-j+1)$. Let $s(j,i)$ and $S(j,i)$, $1\le i\le j$ be the Stirling numbers of the first kind and of the second kind, respectively: $x^{\underline{j}}=\sum\limits_{i=1}^j s(j,i)x^i$ and $x^j=\sum_{i=1}^j S(j,i)x^{\underline{i}}$. Let $r\ge 2$, $n\ge3$ and $x\in[0,1)$. We have
\begin{eqnarray}
Z_n(e_r)(x)&=&\sum_{k=0}^{\infty}\binom{n+k}{k}(1-x)^{n+1}x^k\left(\frac k{n+k}\right)^r\nonumber\\
&=&\sum_{j=1}^r S(r,j)\sum_{k=0}^{\infty}\binom{n+k}{k}(1-x)^{n+1}x^k\frac{k^{\underline{j}}}{(n+k)^r}\nonumber\\
&=&\sum_{j=1}^r S(r,j)\sum_{k=j}^{\infty}\binom{n+k-j}{k-j}(1-x)^{n+1}x^k
\frac{(n+k)^{\underline{j}}}{(n+k)^r}\nonumber\\
&=&\sum_{j=1}^r S(r,j)\sum_{i=1}^j s(j,i)\sum_{k=0}^{\infty}\binom{n+k}{k}(1-x)^{n+1}x^{k+j}\frac{(n+k+j)^i}{(n+k+j)^r}\nonumber\\
&=&\sum_{i=1}^r\sum_{j=i}^rS(r,j)s(j,i)\sum_{k=0}^{\infty}\binom{n+k}{k}(1-x)^{n+1}x^{k+j}\frac1{(n+k+j)^{r-i}}\nonumber\\
&=:&\sum_{i=1}^rT_i.\nonumber
\end{eqnarray}
The permutation of the sums is possible due to the absolute convergence of the involved series. Then
$$T_r=S(r,r)s(r,r)\sum_{k=0}^{\infty}\binom{n+k}{k}(1-x)^{n+1}x^{k+r}=x^r.$$
In order to estimate term $T_{r-1}$ we use the following formula
$$\binom{n+k}{k}\frac1{n+k+j}=\binom{n+k-1}{k}\frac1n+\binom{n+k-2}{k}\left[-\frac j{n(n-1)}+\frac{j(j+1)}{n(n-1)(n+k+j)}\right].$$
Denote
$$\sigma_{n,r}(x)=\sum_{j=r-1}^rS(r,j)s(j,r-1)\sum_{k=0}^{\infty}\binom{n+k-2}{k}(1-x)^{n+1}x^{k+j}\frac{j(j+1)}{n(n-1)(n+k+j)}.$$
We obtain
\begin{eqnarray}
T_{r-1}&=&\sum_{j=r-1}^r S(r,j)s(j,r-1)\left[\frac1n\cdot x^j(1-x)-\frac j{n(n-1)}\cdot x^j(1-x)^2\right]+\sigma_{n,r}(x)\nonumber\\
&=&\binom{r}{2}\sum_{j=r-1}^r(-1)^{r-j-1}x^j(1-x)\left[\frac1n-\frac j{n(n-1)}\cdot (1-x)\right]+\sigma_{n,r}(x)\nonumber\\
&=&\binom{r}{2}x^{r-1}(1-x)^2\left[\frac1n-\frac1{n(n-1)}(r-1-rx)\right]+\sigma_{n,r}(x).\nonumber
\end{eqnarray}

On the other hand we have
\begin{eqnarray}
|\sigma_{n,r}(x)|&\le&\sum_{j=r-1}^r\binom{r}{2}\sum_{k=0}^{\infty}\binom{n+k-2}{k}(1-x)^{n+1}\frac{r(r+1)x^{k+r-1}}{n(n-1)(n+r-1)}\nonumber\\
&=&2\binom{r}{2}\frac{r(r+1)}{n(n-1)(n+r-1)}\cdot x^{r-1}(1-x)^2\nonumber\\
&\le&\frac{r^2(r^2-1)}{n(n^2-1)}\cdot x^{r-1}(1-x)^2.\nonumber
\end{eqnarray}
Finally we have
\begin{eqnarray}
&&\left|\sum_{i=1}^{r-2}T_i\right|=\left|\sum_{i=1}^{r-2}\sum_{j=i}^rS(r,j)s(j,i)\sum_{k=0}^{\infty}\binom{n+k}{k}(1-x)^{n+1}x^{k+j}\cdot\frac1{(n+k+j)^{r-i}}\right|\nonumber\\
&=&\left|\sum_{i=1}^{r-2}\sum_{j=i}^rS(r,j)s(j,i)\frac1{n(n-1)}
\sum_{k=0}^{\infty}\binom{n+k-2}{k}(1-x)^{n+1}x^{k+j}\frac{(n+k)^{\underline{2}}}{(n+k+j)^{r-i}}\right|\nonumber\\
&\le&\frac1{n(n-1)}\sum_{i=1}^{r-2}\sum_{j=i}^r|S(r,j)s(j,i)|x^j(1-x)^2 
\nonumber\\
&\le&\frac{x(1-x)}{n(n-1)}\sum_{i=1}^{r-2}\sum_{j=i}^r|S(r,j)s(j,i)|.\nonumber
\end{eqnarray}
\par From relations above we obtain:
\begin{equation}\label{e33}
Z_n(e_r)(x)=x^r+\frac1n\binom{r}{2}x^{r-1}(1-x)^2+\varepsilon_{n,r}(x)\;
\mbox{ with }\;|\varepsilon_{n,r}(x)|\le\frac{x(1-x)}{n(n-1)}\cdot C_r,
\end{equation}
where $C_r$ is a constant depending only on $r$. This is obviously true also for $x=1$.
\par In the case $r=2$ we obtain an estimate  more refined than (\ref{e33}), which leads to:
\begin{equation}
m_n^2(x)=x(1-x)^2\left[\frac1n-\frac{1-2x}{n(n-1)}+\rho_{n}(x)\right]\;\mbox{with}\;|\rho_n(x)|\le\frac {12}{n(n^2-1)}.
\end{equation}
Denote $\mu_n(x)=\frac12((1-x)\rho_n(x)+x\rho_n(1-x))$. We obtain $M_n^2(x)=\psi(x)\alpha_n(x)$, where
$$\alpha_n(x)=\frac1{2n}-\frac{(1-2x)^2}{2n(n-1)}+\mu_n(x)\;\mbox{ with }\;|\mu_n(x)|\le\frac6{n(n^2-1)}.$$
\par Also, for $t,x\in[0,1]$, we obtain
\begin{eqnarray}
|\alpha_n(t)-\alpha_n(x)|&\le&\frac2{n(n-1)}\cdot|(t-x)(t+x-1)|+|\mu_n(t)|+|\mu_n(x)|\nonumber\\
&\le&\frac{2|t-x|}{n(n-1)}+\frac{12}{n(n^2-1)}.\nonumber
\end{eqnarray}
Note that $Z_n^{\star}(\psi(e_1-xe_0)^2)(x)\le\frac14Z_n^{\star}((e_1-xe_0)^2)(x)\le\frac14\cdot \psi(x)\|\alpha_n\|$ and 
$Z_n^{\star}(\psi)(x)\le\psi(x)$. From the relation above we obtain
\begin{eqnarray}
Z^{\star}_n(\psi|\alpha_n-\alpha_n(x)e_0|)(x)
&\le&\frac2{n(n-1)}\cdot Z_n^{\star}(\psi|e_1-xe_0|)(x)+\frac{12}{n(n^2-1)}\cdot Z_n^{\star}(\psi)(x)\nonumber\\
&\le&\frac2{n(n-1)}\sqrt{Z_n^{\star}(\psi)(x)Z_n^{\star}(\psi(e_1-xe_0)^2)(x)}+\frac{12\psi(x)}{n(n^2-1)}\nonumber\\
&\le&\frac{\psi(x)}{n(n-1)}\cdot\left[\sqrt{\|\alpha_n\|}+\frac{12}{n+1}\right].\nonumber
\end{eqnarray}
Then relation (\ref{e55}) is immediate.
\par From relation (\ref{e33}), for $r=2,3,4$, we obtain, that there is a bounded function $q_n(x)$ with $|q_n(x)|\le M$, $x\in[0,1]$, such that 
$$m_n^4(x)=\frac{x(1-x)}{n(n-1)}\cdot q_n(x).$$
Since $M_n^4(x)=\frac{x(1-x)}{2n(n-1)}\cdot(q_n(x)+q_n(1-x))$, it follows 
$M_n^4(x)\le \frac{M}{n(n-1)\nu_n}\cdot M_n^2(x)$. Therefore, relation~(\ref{e22}) is valid.
\par Consequently, all the conditions in Theorem~\ref{T4} and Theorem~\ref{T5} are satisfied for the sequence of operators $(Z_n^{\star})_{n\ge3}$.

\par\ \\
Ulrich Abel\\
Technische Hochschule Mittelhessen, University of Applied Sciences,\\ 
Department MND, Wilhelm-Leuschner-Stra\ss e 13,\\
61169 Friedberg, Germany\\
e-mail: Ulrich.Abel@mnd.thm.de\\
\par\ \\
Mircea Ivan\\
Technical University of Cluj-Napoca,\\
Department of Mathematics, Str. Memorandumului nr. 28,\\ 
400114 Cluj-Napoca, Romania\\
e-mail: Mircea.Ivan@math.utcluj.ro\\
\par\ \\
Radu P\u alt\u anea\\
"Transilvania" University of Bra\c sov\\
Faculty of Mathematics and Computer Science,\\
Str. Iuliu Maniu nr. 50,\\
500091 Bra\c sov, Romania\\
e-mail: radupaltanea@yahoo.com
\end{document}